\documentclass[12pt]{article}
\usepackage{graphicx}
\usepackage{amsmath}
\usepackage{citesort}
\usepackage{amsfonts}
\usepackage{amssymb}
\begin{document}

\begin{center}
{\textbf{{\Large Explicit inversion formulas for the spherical mean Radon
transform}}}\\[2mm]

\medskip

L. Kunyansky \\[2mm]

University of Arizona, Tucson\\[5mm]
\end{center}

\begin{abstract}
We derive explicit formulas for the reconstruction of a function
from its integrals over a family of spheres, or for the
inversion of the spherical mean Radon transform. Such formulas
are important for problems of thermo- and photo- acoustic tomography.
A closed-form inversion formula of a filtration-backprojection
type is found for the case when the centers of the integration spheres
lie on a sphere in $\mathbb{R}^{n}$ surrounding the support of
the unknown function. An explicit series solution is presented
for the case when the centers of the integration spheres
lie on a general closed surface.
\end{abstract}

\section*{Introduction}

The problem of the reconstruction of a function from its spherical integrals
(or means) have recently attracted attention of researchers due to its
connection to the thermo-acoustic and photo-acoustic tomography
\cite{kruger1,kruger,XuWang,XuWang0}. In these imaging modalities, the object
of interest is illuminated by a short electromagnetic pulse which causes a
fast expansion of the tissue. The intensity of the resulting ultrasound wave
is recorded by a set of detectors surrounding the object. The local intensity
of such expansion is of significant medical interest: it depends on the
physical properties of the tissue (such as, for example, water content) and
its anomaly can be indicative of tumors. Under certain simplifying assumptions
the measurements can be represented by the integrals of the expansion
intensity over the spheres with the centers at the detectors locations. The
reconstruction of the local properties from these integrals is equivalent to
the inversion of the spherical mean Radon transform.

An introduction to the subject can be found in \cite{kruger,kruger1,kuchrev};
for the important results on the injectivity of the spherical Radon transform
and the corresponding range conditions we refer the reader to
\cite{AQ,trinity,Finch,AK1,Finch1}. In the present paper we concentrate on
explicit inversion formulas which are important from both theoretical and
practical points of view. Most of the known formulas of this sort pertain to
the spherical acquisition geometry, i.e. to the situation when centers of the
integration spheres (the positions of the detectors) lie on a sphere
surrounding the body. Such are the series solutions for 2-D and 3-D presented
in \cite{Norton2D,Nort3D,XuWang0}. More desirable backprojection-type formulas
were derived in \cite{Finch} for odd-dimensional spaces, and implemented in
\cite{Haik}. A different explicit formula for the spherical acquisition
geometry valid in 3-D was found in \cite{XuWang} (together with formulas for
certain unbounded acquisition surfaces).

In this paper we present a set of closed-form inversion formulas for the
spherical geometry and a series solution for certain other measuring surfaces.
Our formulas for the spherical case are of the filtration-backprojection type;
they are valid in $\mathbb{R}^{n}$,
$n\geq2.$\\ Such formulas for the
even-dimensional cases were not known previously. (A set of different
inversion formulas for the even-dimensional case was announced by D. Finch
\cite{Finch2} during tomography meeting in Oberwolfach in August, 2006, where
our results where also presented for the first time.)

The spherical case is discussed in Sections~\ref{formulation}
through~\ref{particular}. A series solution for a general acquisition geometry
is presented in Section~\ref{other}.

\section{Formulation of the problem\label{formulation}}

Suppose that $C_{0}^{1}$ function $f(\mathbf{x})$, $\mathbf{x\in}%
\mathbb{R}^{n},\quad n\geq2$ is compactly supported within the closed ball $B$
of radius $R$ centered at the origin. We will denote the boundary of the ball
by $\partial B.$ Our goal is to reconstruct $f(\mathbf{x})$ from its
projections $g(\mathbf{z},r)$ defined as the integrals of $f(\mathbf{x})$ over
the spheres of radius $r$ centered at $\mathbf{z}$:%

\[
g(\mathbf{z},r)=\int\limits_{\mathbb{S}^{n-1}}f(\mathbf{z}+r\hat{t}%
)r^{n-1}ds(\hat{t}\mathbf{).}%
\]
where $\mathbb{S}^{n-1}$ is the unit sphere in $\mathbb{R}^{n},$ $\hat{t}$ is
a unit vector, and $ds$ is the normalized measure in $\mathbb{R}^{n}$.
Projections are assumed to be known for all $\mathbf{z}\in\partial B,$ $0\leq
r\leq2R$ (integrals for $r>2R$ automatically equal zero, since the
corresponding integration spheres do not intersect the support of the
function). In the following two sections we will present an explicit formula
of backprojection type that solves this reconstruction problem.

\section{Derivation}

Our derivation is based on certain properties of the solutions of the
Helmholtz equation in $\mathbb{R}^{n}$%
\[
\Delta h(\mathbf{x})+\lambda^{2}h(\mathbf{x})=0.
\]
For this equation the free space Green's function $\Phi(\mathbf{x,y},\lambda)$
satisfying radiation boundary condition is described (see, for example
\cite{Agmon}) by the formula
\[
\Phi(\mathbf{x,y},\lambda)=\frac{i}{4}\left(  \frac{\lambda}{2\pi
|\mathbf{x}-\mathbf{y}|}\right)  ^{n/2-1}H_{n/2-1}^{(1)}(\lambda
|\mathbf{x}-\mathbf{y}|),
\]
where $H_{n/2-1}^{(1)}(t)$ is the Hankel function of the first kind and of
order $n/2-1.$ To simplify the notation, we introduce functions $J(t),N(t),$
and $H(t),$ defined by the following formulas
\begin{align*}
J(t)  &  =\frac{J_{n/2-1}(t)}{t^{n/2-1}},\\
N(t)  &  =\frac{N_{n/2-1}(t)}{t^{n/2-1}},\\
H(t)  &  =\frac{H_{n/2-1}^{(1)}(t)}{t^{n/2-1}}=J(t)+iN(t),
\end{align*}
where $J_{n/2-1}(t)$ and $N_{n/2-1}(t)$ are respectively the Bessel and
Neumann functions of order $n/2-1.$ In this notation Green's function
$\Phi(\mathbf{x,y},\lambda)$ can be re-written in a simpler form:
\begin{align*}
\Phi(\mathbf{x,y},\lambda)  &  =ic(\lambda,n)H(\lambda|\mathbf{x}%
-\mathbf{y}|)\\
&  =c(\lambda,n)\left[  iJ(\lambda|\mathbf{x}-\mathbf{y}|)-N(\lambda
|\mathbf{x}-\mathbf{y}|)\right]
\end{align*}
where $c(\lambda,n)$ is a constant for a fixed value of $\lambda$:
\[
c(\lambda,n)=\frac{\lambda^{n-2}}{4(2\pi)^{n/2-1}}.
\]
We note that function $J(\lambda|\mathbf{x}|)$ is a solution of the Helmholtz
equation for all $\mathbf{x}\in\mathbb{R}^{n},$ while $N(\lambda|\mathbf{x}|)$
solves this equation in $\mathbb{R}^{n}\backslash\{0\}.$

In other to derive the inversion formula, we utilize the following integral
representation of $f$ in the form of a convolution with $J(\lambda
|\mathbf{y}-\mathbf{x}|)$:
\begin{equation}
f(\mathbf{y})=\frac{1}{(2\pi)^{n/2}}\int\limits_{\mathbb{R}^{+}}\left(
\int\limits_{\mathbb{R}^{n}}f(\mathbf{x})J(\lambda|\mathbf{y}-\mathbf{x}%
|)dx\right)  \lambda^{n-1}d\lambda. \label{intrepr}%
\end{equation}
The above equation easily follows from the Fourier representation of $f$%
\begin{align*}
f(0)  &  =\frac{1}{(2\pi)^{n}}\int\limits_{\mathbb{R}^{n}}\int
\limits_{\mathbb{R}^{n}}e^{-i\mathbf{x}\cdot\xi}f(\mathbf{x})d\mathbf{x}d\xi\\
&  =\frac{1}{(2\pi)^{n}}\int\limits_{\mathbb{R}^{+}}\int\limits_{\mathbb{R}%
^{n}}f(\mathbf{x})\left[  \int\limits_{\mathbb{S}^{n-1}}e^{-i\lambda
\mathbf{x}\cdot\hat{\xi}}d\hat{\xi}\right]  d\mathbf{x}\lambda^{n-1}d\lambda
\end{align*}
and from the well-known integral representation for $J(|\mathbf{u}|)$
\cite{Stein}
\[
J(|\mathbf{u}|)=\frac{1}{(2\pi)^{n/2}}\int\limits_{\mathbb{S}^{n-1}%
}e^{i\mathbf{u}\cdot\hat{\xi}}d\hat{\xi}.
\]

Let us denote the inner integral in (\ref{intrepr}) by $G_{J}(\mathbf{y}%
,\lambda)$
\begin{equation}
G_{J}(\mathbf{y},\lambda)=\int\limits_{\mathbb{R}^{n}}f(\mathbf{x}%
)J(\lambda|\mathbf{y}-\mathbf{x}|)dx. \label{gj}%
\end{equation}
Similarly to the kernel $J(\lambda|\mathbf{y}-\mathbf{x}|)$ of this
convolution, function $G_{J}(\mathbf{y},\lambda)$ is an entire solution of the
Helmholtz equation. Boundary values of this function $G_{J}(\mathbf{z}%
,\lambda),$ $\mathbf{z}\in\partial B$ are easily computable from projections:
\[
G_{J}(\mathbf{z},\lambda)=\int\limits_{B}f(\mathbf{x})J_{0}(\lambda
|\mathbf{z}-\mathbf{x}|)d\mathbf{x}=\int\limits_{0}^{2R}J_{0}(\lambda
r)g(\mathbf{z},r)dr.
\]
If $\lambda$ is not in the spectrum of the Dirichlet Laplacian on $B,$
$G_{J}(\mathbf{y},\lambda)$ is completely determined by its boundary values
and can be found by solving numerically the corresponding Dirichlet problem for
the Helmholtz equation. Then $f(\mathbf{x})$ can be reconstructed from
equation (\ref{intrepr}). Unfortunately, such a solution would not have an
explicit form.

In order to obtain an explicit formula for $G_{J}(\mathbf{y},\lambda)$ we will
utilize a Helmholtz representation for $J(\lambda|\mathbf{y}-\mathbf{x}|);$ it
results from an application of Green's formula and has the following form:
\begin{align*}
J(\lambda|\mathbf{y}-\mathbf{x}|)  &  =\int\limits_{\partial B}\left[
J(\lambda|\mathbf{z}-\mathbf{x}|\frac{\partial}{\partial\mathbf{n}%
_{\mathbf{z}}}\Phi(\mathbf{y,z,}\lambda)\right. \\
&  \left.  -\Phi(\mathbf{y,z},\lambda)\frac{\partial}{\partial\mathbf{n}%
_{\mathbf{z}}}J(\lambda|\mathbf{z}-\mathbf{x}|)\right]  ds(\mathbf{z),}%
\end{align*}
or
\begin{align}
J(\lambda|\mathbf{y}-\mathbf{x}|)  &  =-c(\lambda,n)\int\limits_{\partial
B}\left[  J(\lambda|\mathbf{z}-\mathbf{x}|\frac{\partial}{\partial
\mathbf{n}_{\mathbf{z}}}N(\lambda|\mathbf{y}-\mathbf{z}|)\right. \nonumber\\
&  \left.  -N(\lambda|\mathbf{y}-\mathbf{z}|)\frac{\partial}{\partial
\mathbf{n}_{\mathbf{z}}}J(\lambda|\mathbf{z}-\mathbf{x}|)\right]
ds(\mathbf{z).} \label{helmh}%
\end{align}
Such a representation is valid for any bounded single-connected domain with
sufficiently regular boundary. A straightforward substitution of
equation~(\ref{helmh}) into (\ref{gj}) leads to a boundary value
representation for $G_{J}(\mathbf{y},\lambda)$ involving the normal derivative
the latter function. Unlike the boundary values of $G_{J}(\mathbf{z}%
,\lambda),$ the normal derivative $\frac{\partial}{\partial\mathbf{n}%
_{\mathbf{z}}}G_{J}(\mathbf{z},\lambda)$ cannot be explicitly computed from
projections $g(\mathbf{z},r).$

This difficulty can be circumvented by modifying the Helmholtz representation
as described below. We notice that in the special case of a spherical domain
the second integral in the formula~(\ref{helmh})
\begin{equation}
I(\mathbf{x,y})=\int\limits_{\partial B}N(\lambda|\mathbf{y}-\mathbf{z}%
|)\frac{\partial}{\partial\mathbf{n}_{\mathbf{z}}}J(\lambda|\mathbf{z}%
-\mathbf{x}|)ds(\mathbf{z)} \label{symm}%
\end{equation}
is a symmetric function of its arguments i.e. that
\[
I(\mathbf{x,y})=I(\mathbf{y,x}).
\]
The proof of this fact is presented in the Appendix. Using this symmetry we
obtain a modified Helmholtz representation for $J(\lambda|\mathbf{y}%
-\mathbf{x}|)$:
\begin{align}
J(\lambda|\mathbf{y}-\mathbf{x}|)  &  =-c(\lambda,n)\int\limits_{\partial
B}\left[  J(\lambda|\mathbf{z}-\mathbf{x}|\frac{\partial}{\partial
\mathbf{n}_{\mathbf{z}}}N(\lambda|\mathbf{y}-\mathbf{z}|)\right. \nonumber\\
&  \left.  -N(\lambda|\mathbf{z}-\mathbf{x}|)\frac{\partial}{\partial
\mathbf{n}_{\mathbf{z}}}J(\lambda|\mathbf{y}-\mathbf{z}|)\right]
ds(\mathbf{z).} \label{modhelm}%
\end{align}
\bigskip Now the substitution of equation~(\ref{modhelm}) into (\ref{gj})
yields
\begin{align*}
\int\limits_{B}f\mathbf{(x)}J(\lambda|\mathbf{y}-\mathbf{x}|)dx  &
=-c(\lambda,n)\int\limits_{\partial B}\left[  \left(  \int\limits_{B}%
f\mathbf{(x)}J(\lambda|\mathbf{z}-\mathbf{x}|)dx\right)  \frac{\partial
}{\partial\mathbf{n}_{\mathbf{z}}}N(\lambda|\mathbf{y}-\mathbf{z}|)\right. \\
&  \left.  -\left(  \int\limits_{B}f\mathbf{(x)}N(\lambda|\mathbf{z}%
-\mathbf{x}|)dx\right)  \frac{\partial}{\partial\mathbf{n}_{\mathbf{z}}%
}J(\lambda|\mathbf{y}-\mathbf{z}|)\right]  ds(\mathbf{z),}%
\end{align*}
where the inner integrals are easily computable from the projections
$g(\mathbf{z},t)$:
\begin{align*}
\int\limits_{B}f\mathbf{(x)}J(\lambda|\mathbf{z}-\mathbf{x}|)dx  &
=\int\limits_{0}^{2R}J(\lambda t)g(\mathbf{z},t)dt,\\
\int\limits_{B}f\mathbf{(x)}N(\lambda|\mathbf{z}-\mathbf{x}|)dx  &
=\int\limits_{0}^{2R}N(\lambda t)g(\mathbf{z},t)dt.
\end{align*}
Thus, the convolution of $f$ and $J$ can be reconstructed from the projections
as follows:
\begin{align}
\int\limits_{B}f\mathbf{(x)}J(\lambda|\mathbf{y}-\mathbf{x}|)dx  &
=-c(\lambda,n)\int\limits_{\partial B}\left[  \left(  \int\limits_{0}%
^{2R}J(\lambda t)g(\mathbf{z},t)dt\right)  \frac{\partial}{\partial
\mathbf{n}_{\mathbf{z}}}N(\lambda|\mathbf{y}-\mathbf{z}|)\right. \nonumber\\
&  \left.  -\left(  \int\limits_{0}^{2R}N(\lambda t)g(\mathbf{z},t)dt\right)
\frac{\partial}{\partial\mathbf{n}_{\mathbf{z}}}J(\lambda|\mathbf{y}%
-\mathbf{z}|)\right]  ds(\mathbf{z)}\nonumber\\
&  =c(\lambda,n)\mathrm{div}\int\limits_{\partial B}\mathbf{n}(\mathbf{z}%
)\left[  \left(  \int\limits_{0}^{2R}J(\lambda t)g(\mathbf{z},t)dt\right)
N(\lambda|\mathbf{y}-\mathbf{z}|)\right. \nonumber\\
&  \left.  -\left(  \int\limits_{0}^{2R}N(\lambda t)g(\mathbf{z},t)dt\right)
J(\lambda|\mathbf{y}-\mathbf{z}|)\right]  ds(\mathbf{z).} \label{convol}%
\end{align}
Finally, by combining equations~(\ref{intrepr}) and (\ref{convol}), one
arrives at the following inversion formula
\[
f(\mathbf{y})=\frac{1}{4(2\pi)^{n-1}}\mathrm{div}\int\limits_{\partial
B}\mathbf{n}(\mathbf{z})h(\mathbf{z},|\mathbf{y}-\mathbf{z}|)ds(\mathbf{z),}%
\]
where
\begin{align*}
h(\mathbf{z},t)  &  =\int\limits_{\mathbb{R}^{+}}\left[  N(\lambda t)\left(
\int\limits_{0}^{2R}J(\lambda t^{\prime})g(\mathbf{z},t^{\prime})dt^{\prime
}\right)  \right. \\
&  -J(\lambda t)\left.  \left(  \int\limits_{0}^{2R}N(\lambda t^{\prime
})g(\mathbf{z},t^{\prime})dt^{\prime}\right)  \right]  \lambda^{2n-3}d\lambda.
\end{align*}

\section{Particular cases\label{particular}}

\subsection{2-D case}

From the point of view of practical applications the two- and three-
dimensional cases are the most important ones. In 2-D, $J(t)=J_{0}(t),$
$N(t)=N_{0}(t),$ and the inversion formula has the following form
\begin{equation}
f(\mathbf{y})=\frac{1}{8\pi}\mathrm{div}\int\limits_{\partial B}%
\mathbf{n}(\mathbf{z})h(\mathbf{z},|\mathbf{y}-\mathbf{z}|)dl(\mathbf{z)},
\label{backpro}%
\end{equation}
where
\begin{align}
h(\mathbf{z},t)  &  =\int\limits_{\mathbb{R}^{+}}\left[  N_{0}(\lambda
t)\left(  \int\limits_{0}^{2R}J_{0}(\lambda t^{\prime})g(\mathbf{z},t^{\prime
})dt^{\prime}\right)  \right. \nonumber\\
&  -\left.  J_{0}(\lambda t)\left(  \int\limits_{0}^{2R}N_{0}(\lambda
t^{\prime})g(\mathbf{z},t^{\prime})dt^{\prime}\right)  \right]  \lambda
d\lambda. \label{filtration}%
\end{align}
Equation (\ref{backpro}) is a backprojection followed by the divergence
operator. It is worth noticing that such a divergence form of a reconstruction
formula is not unusual; it also naturally occurs in reconstruction formulas
for the attenuated Radon transform (see, for instance \cite{Novikov}).

Equation (\ref{filtration}) represents the filtration step of the algorithm.
In order to better understand the nature of this operator we re-write
(\ref{filtration}) in the form
\begin{align}
h(\mathbf{z},t)  &  =-\frac{i}{2}\int\limits_{\mathbb{R}^{+}}\left[
H_{0}^{(1)}(\lambda t)\left(  \int\limits_{0}^{2R}\overline{H_{0}%
^{(1)}(\lambda t^{\prime})}g(\mathbf{z},t^{\prime})dt^{\prime}\right)  \right.
\nonumber\\
&  -\left.  \overline{H_{0}^{(1)}(\lambda t)}\left(  \int\limits_{0}^{2R}%
H_{0}^{(1)}(\lambda t^{\prime})g(\mathbf{z},t^{\prime})dt^{\prime}\right)
\right]  \lambda d\lambda. \label{hankfilt}%
\end{align}
and recall that for large values of the argument $t$ the Hankel function
$H_{0}^{(1)}(t)$ has the following asymptotic expansion (\cite{Watson})
\[
H_{0}^{(1)}(t)=\left(  \frac{2}{\pi t}\right)  ^{\frac{1}{2}}e^{-\frac{i}%
{2}\pi}e^{it}.
\]
In the situation when the support of the function $f(\mathbf{x})$ remains
bounded and the radius $R$ of the ball $B$ becomes large, the inner and outer
integrals in (\ref{hankfilt}) reduce to the direct and inverse Fourier
transforms, with one of the terms corresponding to the positive frequencies
and the other to the negative ones. Thus, since the two terms have opposite
signs, in the asymptotic limit of large $R$ operator (\ref{hankfilt}) equals
(up to a constant factor) to the Hilbert transform of $g(\mathbf{z},\cdot)$.

\subsection{3-D case}

In the three-dimensional case
\begin{align*}
J(t)  &  =\frac{J_{1/2}(\lambda r)}{\sqrt{\lambda r}}=\sqrt{\frac{2}{\pi}%
}j_{0}(\lambda r),\\
N(t)  &  =\frac{N_{1/2}(\lambda r)}{\sqrt{\lambda r}}=\sqrt{\frac{2}{\pi}%
}n_{0}(\lambda r),
\end{align*}
where $j_{0}(t)$ and $n_{0}(t)$ are the spherical Bessel and Neumann functions
respectively. The formula takes the following form
\begin{equation}
f\mathbf{(y)=}\frac{1}{16\pi^{2}}\mathrm{div}\int\limits_{|z|=R}%
\mathbf{n(z)}h\mathbf{(z,|y-z|)}ds\mathbf{(z),} \label{form3d}%
\end{equation}
with
\begin{align*}
h(\mathbf{z},t)  &  =\frac{2}{\pi}\int\limits_{\mathbb{R}^{+}}\left[
n_{0}(\lambda t)\left(  \int\limits_{0}^{2R}j_{0}(\lambda t^{\prime
})g(\mathbf{z},t^{\prime})dt^{\prime}\right)  \right. \\
&  -\left.  j_{0}(\lambda t)\left(  \int\limits_{0}^{2R}n_{0}(\lambda
t^{\prime})g(\mathbf{z},t^{\prime})dt^{\prime}\right)  \right]  \lambda
^{3}d\lambda.
\end{align*}
In this case, however, a further simplification is possible since $j_{0}(t)$
and $n_{0}(t)$ have a simple representation in terms of trigonometric
functions:
\[
j_{0}(t)=\frac{\sin t}{t},\quad n_{0}(t)=-\frac{\cos t}{t}.
\]
The substitution of these trigonometric expressions into the inversion formula
leads to a significantly simpler formula:
\begin{align}
h(\mathbf{z},t)  &  =-\frac{2}{\pi t}\int\limits_{\mathbb{R}^{+}}\cos(\lambda
t)\left[  \int\limits_{0}^{2R}\sin(\lambda t^{\prime})\frac{g(\mathbf{z}%
,t^{\prime})}{t^{\prime}}dt^{\prime}\right]  \lambda d\lambda\nonumber\\
&  +\frac{2}{\pi t}\int\limits_{\mathbb{R}^{+}}\sin(\lambda t)\left[
\int\limits_{0}^{2R}\cos(\lambda t^{\prime})\frac{g(\mathbf{z},t^{\prime}%
)}{t^{\prime}}dt^{\prime}\right]  \lambda d\lambda\nonumber\\
&  =-\frac{2}{\pi t}\frac{d}{dt}\int\limits_{\mathbb{R}^{+}}\sin(\lambda
t)\left[  \int\limits_{0}^{2R}\sin(\lambda t^{\prime})\frac{g(\mathbf{z}%
,t^{\prime})}{t^{\prime}}dt^{\prime}\right]  d\lambda\nonumber\\
&  -\frac{2}{\pi t}\frac{d}{dt}\int\limits_{\mathbb{R}^{+}}\cos(\lambda
t)\left[  \int\limits_{0}^{2R}\cos(\lambda t^{\prime})\frac{g(\mathbf{z}%
,t^{\prime})}{t^{\prime}}dt^{\prime}\right]  d\lambda\nonumber\\
&  =-\frac{2}{t}\frac{d}{dt}\frac{g(\mathbf{z},t)}{t}, \label{deriv}%
\end{align}
where we took into account the fact that the Fourier sine and cosine
transforms are self-invertible. By combining (\ref{deriv}) and (\ref{form3d})
our inversion formula can be re-written in the form
\[
f(\mathbf{y})=-\frac{1}{8\pi^{2}}\mathrm{div}\int\limits_{\partial
B}\mathbf{n}(\mathbf{z})\left( \frac{1}{t} \frac{d}{dt}\frac{g(\mathbf{z},t)}{t}\right)  \left.
{\phantom{\rule{1pt}{8mm}}}\right|  _{t=|\mathbf{z}-\mathbf{y}|}%
ds(\mathbf{z}).
\]
This expression is equivalent to one of the formulas derived in \cite{XuWang}
for the 3-D case.

\section{Inversion of the spherical mean Radon transform in other
geometries\label{other}}

In such applications as photo- and thermo- acoustic tomography, the designer
of the measuring system has a freedom in selecting the detectors' locations.
The detectors (the centers of the integration spheres) do not have to lie on a
sphere. In this section we present reconstruction formulas for the case when
the measuring surface is a boundary of certain other domains. Namely, our
method works for the domains whose eigenfunctions $u_{m}(\mathbf{x})$ of the
(zero) Dirichlet Laplacian are explicitly known. Such are, for example, the
domains for which the eigenfunctions can be found by the separation of
variables, i.e. sphere, annulus, cube, and certain subsets of those, and
crystallographic domains (see \cite{berard1,berard2}).

The proof of the range theorem in~\cite{trinity} involves implicitly a
reconstruction procedure also based on eigenfunction expansions. Unlike the
present method, that procedure would involve division of analytic functions
that have countable number of zeros. While the range theorem guarantees
cancellation of these zeros when the data are in the range of the direct
transform, a stable numerical implementation of such division would be
complicated if not impossible. (Similarly, the series solution
of~\cite{Norton2D} for 2-D circular geometry involves division by Bessel
functions.) The technique we present below does not require such divisions.

Suppose $\lambda_{m}^{2}$, $u_{m}(\mathbf{x})$ are the eigenvalues and
eigenfunctions of the (negative) Dirichlet Laplacian on a bounded domain
$\Omega$ with zero boundary conditions, i.e.
\begin{align*}
\Delta u_{m}(\mathbf{x})+\lambda_{m}^{2}u_{m}(\mathbf{x}) &  =0,\qquad
\mathbf{x}\in\Omega,\quad\Omega\subseteq\mathbb{R}^{n},\\
u_{m}(\mathbf{x}) &  =0,\qquad\mathbf{x}\in\partial\Omega.
\end{align*}
As in the previous sections, we would like to reconstruct a function
$f(\mathbf{x})\in$ $L^{2}(\Omega)$ from the known values of its spherical
integrals $g(\mathbf{z},r)$ with the centers on $\partial\Omega$:
\[
g(\mathbf{z},r)=\int\limits_{\mathbb{S}^{n-1}}f(\mathbf{z}+r\mathbf{\hat{s}}%
)r^{n-1}d\mathbf{\hat{s}},\qquad\mathbf{z}\in\partial\Omega.
\]
We notice that $u_{m}(\mathbf{x})$ is a solution of the Dirichlet problem for
the Helmholtz equation and thus admits the Helmholtz representation
\begin{align}
u_{m}(\mathbf{x}) &  =\int_{\partial\Omega}\Phi(\mathbf{x},\mathbf{z}%
,\lambda_{m})\frac{\partial}{\partial\mathbf{n}}u_{m}(\mathbf{z}%
)ds(\mathbf{z)}\nonumber\\
&  =ic(\lambda_{m},n)\int_{\partial\Omega}H(\lambda_{m}|\mathbf{x}%
-\mathbf{z}|)\frac{\partial}{\partial\mathbf{n}}u_{m}(\mathbf{z}%
)ds(\mathbf{z)}\qquad\mathbf{x}\in\Omega.\label{helmdiscr}%
\end{align}

On the other hand, eigenfunctions $\left\{  u_{m}(\mathbf{x})\right\}
_{0}^{\infty}$ form an orthonormal basis in $L_{2}(\Omega).$ Therefore
$f(\mathbf{x})$ can be represented (in $L^{2}$sense) by the series
\begin{equation}
f(\mathbf{x})=\sum_{m=0}^{\infty}\alpha_{m}u_{m}(\mathbf{x})
\label{fourierser}%
\end{equation}
with
\begin{equation}
\alpha_{m}=\int_{\Omega}u_{m}(\mathbf{x})f(\mathbf{x})d\mathbf{x.}
\label{serkoef}%
\end{equation}
The reconstruction formula will result if we substitute
representation~(\ref{helmdiscr}) into (\ref{serkoef}) and change the order of
integrations
\begin{align}
\alpha_{m}  &  =\int_{\Omega}u_{m}(\mathbf{x})f(\mathbf{x})d\mathbf{x}%
\nonumber\\
&  =ic(\lambda_{m},n)\int_{\Omega}\left(  \int_{\partial\Omega}H(\lambda
_{m}|\mathbf{x}-\mathbf{z}|)\frac{\partial}{\partial\mathbf{n}}u_{m}%
(\mathbf{z})ds(\mathbf{z)}\right)  f(\mathbf{x})d\mathbf{x}\nonumber\\
&  =ic(\lambda_{m},n)\int_{\partial\Omega}\left(  \int_{\Omega}H(\lambda
_{m}|\mathbf{x}-\mathbf{z}|)f(\mathbf{x})d\mathbf{x}\right)  \frac{\partial
}{\partial\mathbf{n}}u_{m}(\mathbf{z})ds(\mathbf{z).} \label{serkoef1}%
\end{align}
The change of the integration order is justified by the fact that
eigenfunctions $u_{m}(\mathbf{x})$ are continuous. The inner integral
in~(\ref{serkoef1}) is easily computed from projections
\[
\int_{\Omega}H(\lambda_{m}|\mathbf{x}-\mathbf{z}|)f(\mathbf{x})d\mathbf{x}%
=\int\limits_{\mathbb{R}^{+}}g(\mathbf{z},r)H(\lambda_{m}r)dr,
\]
so that
\begin{equation}
\alpha_{m}=ic(\lambda_{m},n)\int_{\partial\Omega}\left(  \int
\limits_{\mathbb{R}^{+}}g(\mathbf{z},r)H(\lambda_{m}r)dr\right)
\frac{\partial}{\partial\mathbf{n}}u_{m}(\mathbf{z})ds(\mathbf{z).}
\label{serkoef2}%
\end{equation}
With Fourier coefficients $\alpha_{m}$ now known $f(\mathbf{x})$is
reconstructed by summing series~(\ref{fourierser}).

If desired, this solution can be re-written in the form of a
backprojection-type formula:
\begin{align}
f(\mathbf{x})  &  =\sum_{m=0}^{\infty}\alpha_{m}u_{m}(\mathbf{x}%
)=\int_{\partial\Omega}\left(  \sum_{m=0}^{\infty}\alpha_{m}c(\lambda
_{m},n)iH(\lambda_{m}|\mathbf{x}-\mathbf{z}|)\frac{\partial}{\partial
\mathbf{n}}u_{m}(\mathbf{z})\right)  ds(\mathbf{z)}\nonumber\\
&  =\int_{\partial\Omega}h(\mathbf{z},|\mathbf{x}-\mathbf{z|})ds(\mathbf{z),}
\label{serbackpr}%
\end{align}
where
\begin{equation}
h(\mathbf{z},t)=i\sum_{m=0}^{\infty}c(\lambda_{m},n)\alpha_{m}H(\lambda
_{m}t)\frac{\partial}{\partial\mathbf{n}}u_{m}(\mathbf{z}), \label{serfiltr}%
\end{equation}
and coefficients $a_{m}$ are computed using equation~(\ref{serkoef2}). In the
above formula equation~(\ref{serbackpr}) is clearly a backprojection operator,
and~(\ref{serfiltr}) is a filtration. However, the latter operator is now
represented by a series rather than by a closed form expression.

The series solution described above have an interesting property not possessed
(to the best of our knowledge) by any other explicit reconstruction technique.
Let us consider a slightly more general problem. Suppose that region $\Omega$
is a proper subset of a larger region $\Omega_{1}$ ($\Omega\subset\Omega_{1}$)
and that a $L^{2}$ function $F$ is defined on $\Omega_{1}.$ We will denote the
restriction of $F$ on $\Omega$ by $f,$ i.e.
\[
f(\mathbf{x})=\left\{
\begin{array}
[c]{ccc}%
F(\mathbf{x}) & , &\mathbf{x}\in\Omega\\
0 & , & \mathbf{x}\in\mathbb{R}^{n}\backslash\Omega
\end{array}
\right.  .
\]
We would like to reconstruct $f(\mathbf{x})$ from the integrals $g(\mathbf{z}%
,r)$ of $F$ over spheres with the centers on $\partial\Omega$:
\[
g(\mathbf{z},r)=\int\limits_{\mathbb{S}^{n-1}}F(\mathbf{z}+r\mathbf{\hat{s}%
})r^{n-1}d\mathbf{\hat{s}},\qquad\mathbf{z}\in\partial\Omega.
\]
The difference with the previously considered problem is in that the centers
of the integration spheres are know lying on a surface which is inside the
support $\Omega_{1}$ of the function $F.$ While we are still trying to
reconstruct the restriction $f$ of $F$ to $\Omega$, the integrals we know are
those of $F$ and not of $f.$

It turns out that the solution to this problem is still given by
formulas~(\ref{fourierser}) and~(\ref{serkoef2}) (or equivalently
by~(\ref{serbackpr}), (\ref{serfiltr}), and~(\ref{serkoef2})). Indeed, if we
extend functions $u_{m}(\mathbf{x})$ by $0$ to $\mathbb{R}^{n}\backslash
\Omega$, formula~(\ref{helmdiscr}) holds for all $\mathbf{x} \in \mathbb{R}^n$,
and~(\ref{fourierser})  remains unchanged. In
formula~(\ref{serkoef1}) $f$ can be replaced by $F$ as follows:
\begin{align*}
\alpha_{m}  &  =\int_{\Omega}u_{m}(\mathbf{x})f(\mathbf{x})d\mathbf{x=}%
\int_{\Omega}u_{m}(\mathbf{x})F(\mathbf{x})d\mathbf{x}\\
&  =ic(\lambda_{m},n)\int_{\partial\Omega}\left(  \int_{\Omega}H(\lambda
_{m}|\mathbf{x}-\mathbf{z}|)F(\mathbf{x})d\mathbf{x}\right)  \frac{\partial
}{\partial\mathbf{n}}u_{m}(\mathbf{z})ds(\mathbf{z).}%
\end{align*}
and the inner integral can be computed from projections as before:
\[
\int_{\Omega}H(\lambda_{m}|\mathbf{x}-\mathbf{z}|)F(\mathbf{x})d\mathbf{x}%
=\int\limits_{\mathbb{R}^{+}}g(\mathbf{z},r)H(\lambda_{m}r)dr.
\]
By combining the two above equations we again arrive at the
formula~(\ref{serkoef2}).

To summarize, if the centers of the integration spheres lie on a closed
surface inside of the support of the function, the values of the function
corresponding to the interior of that surface can be stably reconstructed by a
formula containing the normal derivatives of the eigenfunctions of the
Dirichlet Laplacian. If these eigenfunctions are given by an explicit formula,
our technique yields an explicit series solution to this problem.

\section*{Appendix}

In this section we prove that for arbitrary $\mathbf{x},\mathbf{y}\in B$
function $I(\mathbf{x},\mathbf{y})$ defined by equation~(\ref{symm}) is a
symmetric function of its arguments, i.e. that $I(\mathbf{x},\mathbf{y}%
)=I(\mathbf{y},\mathbf{x})$.

First we need to find a closed form expression for a single layer potential of
a spherical shell with the density equal to a spherical harmonic $Y_{l}%
^{(k)}(\hat{x}).$ To this end we introduce functions $u_{\lambda}%
^{k,l}(\mathbf{x})$ and $v_{\lambda}^{k,l}(\mathbf{x})$ defined as follows
\begin{align*}
u_{\lambda}^{k,l}(\mathbf{x})  &  =Y_{l}^{(k)}(\hat{x})J_{(k)}(\lambda
|\mathbf{x}|),\\
v_{\lambda}^{k,l}(\mathbf{x})  &  =Y_{l}^{(k)}(\hat{x})H_{(k)}(\lambda
|\mathbf{x}|),\qquad\hat{x}=\mathbf{x}/|\mathbf{x}|,
\end{align*}
where we, as before, utilize the notation
\begin{align*}
J_{(k)}(t)  &  =\frac{J_{n/2+k-1}(t)}{t^{n/2-1}},\\
N_{(k)}(t)  &  =\frac{N_{n/2+k-1}(t)}{t^{n/2-1}},\\
H_{(k)}(t)  &  =\frac{H_{n/2+k-1}^{(1)}(t)}{t^{n/2-1}}.
\end{align*}
Functions $u_{\lambda}^{k,l}(\mathbf{x})$ are (\cite{Volchkov}) entire
solutions of the Helmholtz equation for all $\mathbf{x}\in\mathbb{R}^{n},$
while $v_{\lambda}^{k,l}(\mathbf{x})$ are radiating solutions of this equation
in $\mathbb{R}^{n}\backslash\{0\}.$ If we apply the Green's theorem in an
exterior of a sphere of radius $r_{0}$ to functions $u_{\lambda}%
^{k,l}(\mathbf{x})$ and $v_{\lambda}^{k,l}(\mathbf{x}),$ we will obtain (for
$|\mathbf{x}|>r_{0})$%
\begin{align*}
0  &  =\int\limits_{|\mathbf{z}|=r_{0}}\left[  u_{\lambda}^{k,l}%
(\mathbf{z})\frac{\partial}{\partial\mathbf{n}_{\mathbf{z}}}\Phi
(\mathbf{x},\mathbf{z},\lambda)-\Phi(\mathbf{x},\mathbf{z},\lambda
)\frac{\partial}{\partial\mathbf{n}_{\mathbf{z}}}u_{\lambda}^{k,l}%
(\mathbf{z})\right]  d\mathbf{z}\\
&  =\int\limits_{|\mathbf{z}|=r_{0}}Y_{l}^{(k)}(\hat{z})\left[  J_{(k)}%
(\lambda r_{0})\frac{\partial}{\partial\mathbf{n}_{\mathbf{z}}}\Phi
(\mathbf{x},\mathbf{z},\lambda)-\Phi(\mathbf{x},\mathbf{z},\lambda)\lambda
J_{(k)}^{\prime}(\lambda r_{0})\right]  d\mathbf{z},
\end{align*}
and
\begin{align*}
v_{\lambda}^{k,l}(\mathbf{x})  &  =\int\limits_{|\mathbf{z}|=r_{0}}\left[
v_{\lambda}^{k,l}(\mathbf{z})\frac{\partial}{\partial\mathbf{n}_{\mathbf{z}}%
}\Phi(\mathbf{x},\mathbf{z},\lambda)-\Phi(\mathbf{x},\mathbf{z},\lambda
)\frac{\partial}{\partial\mathbf{n}_{\mathbf{z}}}v_{\lambda}^{k,l}%
(\mathbf{z})\right]  d\mathbf{z}\\
&  =\int\limits_{|\mathbf{z}|=r_{0}}Y_{l}^{(k)}(\hat{z})\left[  H_{(k)}%
(\lambda r_{0})\frac{\partial}{\partial\mathbf{n}_{\mathbf{z}}}\Phi
(\mathbf{x},\mathbf{z},\lambda)-\Phi(\mathbf{x},\mathbf{z},\lambda)\lambda
H_{(k)}^{\prime}(\lambda r_{0})\right]  d\mathbf{z}%
\end{align*}
where $\hat{z}=\mathbf{z}/|\mathbf{z}|$. By combining the above two equations
one can eliminate the terms with $\frac{\partial}{\partial\mathbf{n}%
_{\mathbf{z}}}$:
\begin{align}
J_{(k)}(\lambda r_{0})v_{\lambda}^{k,l}(\mathbf{x})  &  =\lambda\left[
H_{(k)}(\lambda r_{0})J_{(k)}^{\prime}(\lambda r_{0})-J_{(k)}(\lambda
r_{0})H_{(k)}^{\prime}(\lambda r_{0})\right] \nonumber\\
&  \times\int\limits_{|\mathbf{z}|=r_{0}}Y_{l}^{(k)}(\hat{z})\Phi
(\mathbf{x},\mathbf{z},\lambda)d\mathbf{z.} \label{formula1}%
\end{align}
The expression in the brackets can be simplified using the well known formula
(\cite{Watson}) for the Wronskian of $H_{\alpha}^{(0)}(t)$ and $J_{\alpha}%
(t)$:
\[
H_{\alpha}^{(0)}(t)J_{\alpha}^{\prime}(t)-J_{\alpha}(t)\left(  H_{\alpha
}^{(0)}\right)  ^{\prime}(t)=-\frac{i}{2\pi t}.
\]
Formula (\ref{formula1}) then takes form
\[
J_{(k)}(\lambda r_{0})v_{\lambda}^{k,l}(\mathbf{x})=-\frac{i}{2\pi
r_{0}(\lambda r_{0})^{n-2}}\int\limits_{|\mathbf{z}|=r_{0}}Y_{l}^{(k)}(\hat
{z})\Phi(\mathbf{x},\mathbf{z},\lambda)d\mathbf{z}%
\]
so that the single layer potential we consider is described by the following
equation
\begin{equation}
\int\limits_{|\mathbf{z}|=r_{0}}Y_{l}^{(k)}(\hat{z})\Phi(\mathbf{x}%
,\mathbf{z},\lambda)d\mathbf{z}=i2\pi r_{0}(\lambda r_{0})^{n-2}%
J_{(k)}(\lambda r_{0})H_{(k)}(\lambda|\mathbf{x}|)Y_{l}^{(k)}(\hat{x}).
\label{formula2}%
\end{equation}
By substituting into (\ref{formula2}) the expression for the Green's function
$\Phi(\mathbf{x},\mathbf{z},\lambda)$ in the form
\[
\Phi(\mathbf{x},\mathbf{z},\lambda)=ic(\lambda,n)H(\lambda|\mathbf{x}%
-\mathbf{z}|)
\]
one obtains
\[
\int\limits_{|\mathbf{z}|=r_{0}}Y_{l}^{(k)}(\hat{z})H(\lambda|\mathbf{x}%
-\mathbf{z}|)d\mathbf{z}=\frac{2\pi r_{0}(\lambda r_{0})^{n-2}}{c(\lambda
,n)}J_{(k)}(\lambda r_{0})H_{(k)}(\lambda|\mathbf{x}|)Y_{l}^{(k)}(\hat{x}).
\]
Finally, by separating the real and imaginary parts of the above equation we
arrive at the following two formulas:
\begin{equation}
\int\limits_{|\mathbf{z}|=r_{0}}Y_{l}^{(k)}(\hat{z})J(\lambda|\mathbf{x}%
-\mathbf{z}|)d\mathbf{z}=\frac{2\pi r_{0}(\lambda r_{0})^{n-2}}{c(\lambda
,n)}J_{(k)}(\lambda r_{0})J_{(k)}(\lambda|\mathbf{x}|)Y_{l}^{(k)}(\hat{x}),
\label{jequ}%
\end{equation}%
\begin{equation}
\int\limits_{|\mathbf{z}|=r_{0}}Y_{l}^{(k)}(\hat{z})N(\lambda|\mathbf{x}%
-\mathbf{z}|)d\mathbf{z}=\frac{2\pi r_{0}(\lambda r_{0})^{n-2}}{c(\lambda
,n)}J_{(k)}(\lambda r_{0})N_{(k)}(\lambda|\mathbf{x}|)Y_{l}^{(k)}(\hat{x}),
\label{nequ}%
\end{equation}
valid for $|\mathbf{x}|>r_{0}.$ This completes the preparation for the proof
of the symmetry of $I(\mathbf{x},\mathbf{y}).$

Let us consider function $I(\alpha\hat{x},\beta\hat{y}).$ For fixed values of
$\alpha$ and $\beta$ this is an infinitely smooth function of $\hat{x}$ and
$\hat{y}$ defined on $\mathbb{S}^{n-1}\times\mathbb{S}^{n-1}.$ Consider the
Fourier expansion of $I(\alpha\hat{x},\beta\hat{y})$ in the spherical
harmonics in both variables $\hat{x}$ and $\hat{y}$. Coefficients of such
series are given by the formula
\begin{equation}
a_{k^{\prime,}l^{\prime}}^{k,l}(\alpha,\beta)=\int\limits_{\mathbb{S}^{n-1}%
}\int\limits_{\mathbb{S}^{n-1}}Y_{l}^{k}(\hat{x})Y_{l^{\prime}}^{k^{\prime}%
}(\hat{y})I(\alpha\hat{x},\beta\hat{y})d\hat{x}d\hat{y} \label{ycoefs}%
\end{equation}
with $k=0,1,2,...,\quad0\leq l\leq d_{k},$%
\begin{align*}
d_{k}  &  =%
\genfrac{(}{)}{}{}{n+k-1}{k}%
-%
\genfrac{(}{)}{}{}{n+k-3}{k-2}%
,k\geq2,\\
d_{1}  &  =n,d_{0}=1.
\end{align*}
In order to prove the symmetry $I(\mathbf{x},\mathbf{y})=I(\mathbf{y,x})$ it
is enough to prove that
\begin{equation}
a_{k^{\prime,}l^{\prime}}^{k,l}(\alpha,\beta)=a_{k^{,}l}^{k^{\prime}%
,l^{\prime}}(\beta,\alpha) \label{symcoef}%
\end{equation}
for all relevant values of $k,k^{\prime},l,l^{\prime}.$ By substituting the
expression for $I(\alpha\hat{x},\beta\hat{y})$ into (\ref{ycoefs}) one obtains%

\begin{align*}
a_{k^{\prime,}l^{\prime}}^{k,l}(\alpha,\beta) &  =\int\limits_{\mathbb{S}%
^{n-1}}\int\limits_{\mathbb{S}^{n-1}}Y_{l}^{k}(\hat{x})Y_{l^{\prime}%
}^{k^{\prime}}(\hat{y})\left[  \int\limits_{|\mathbf{z}|=R}N(\lambda
|\mathbf{z}-\alpha\hat{x}|)\frac{\partial}{\partial\mathbf{n}_{\mathbf{z}}%
}J(\beta\hat{y}-\mathbf{z}|)d\mathbf{z}\right]  d\hat{x}d\hat{y}\\
&  =\int\limits_{\mathbb{S}^{n-1}}Y_{l^{\prime}}^{k^{\prime}}(\hat{y})\left[
\int\limits_{|\mathbf{z}|=R}\left(  \int\limits_{\mathbb{S}^{n-1}}Y_{l}%
^{k}(\hat{x})N(\lambda|\mathbf{z}-\alpha\hat{x}|)d\hat{x}\right)
\frac{\partial}{\partial\mathbf{n}_{\mathbf{z}}}J(\lambda|\beta
\hat{y}-\mathbf{z}|d\mathbf{z}\right]  d\hat{y}\\
&  =\frac{1}{\alpha^{n-1}}\int\limits_{\mathbb{S}^{n-1}}Y_{l^{\prime}%
}^{k^{\prime}}(\hat{y})\left[  \int\limits_{|\mathbf{z}|=R}\left(
\int\limits_{|u|=|\mathbf{x}|}Y_{l}^{k}(u/|u|)N(\lambda|\mathbf{z}%
-u|)du\right)  \frac{\partial}{\partial\mathbf{n}_{\mathbf{z}}}J(\lambda
|\beta\hat{y}-\mathbf{z}|d\mathbf{z}\right]  d\hat{y}.
\end{align*}
Utilizing formulas (\ref{jequ}), (\ref{nequ}) and the orthonormality of the
spherical harmonics we find that $a_{k^{\prime,}l^{\prime}}^{k,l}(\alpha
,\beta)=0$ if $k\neq k^{\prime}$ or $l\neq l^{\prime}.$ Otherwise
\[
a_{k^{,}l}^{k,l}(\alpha,\beta)=\lambda\left(  \frac{2\pi\lambda^{n-2}%
}{c(\lambda,n)}\right)  ^{2}R^{n-1}J_{(k)}(\lambda\alpha)J_{(k)}(\lambda
\beta)N_{(k)}(\lambda R)J_{(k)}^{\prime}(\lambda R).
\]
Inspection of the above formula shows that coefficients $a_{k^{\prime
,}l^{\prime}}^{k,l}(\alpha,\beta)$ indeed satisfy (\ref{symcoef}) and, thus,
that $I(\mathbf{x},\mathbf{y})=I(\mathbf{y},\mathbf{x})$.

\section{Acknowledgments}

The author would like to thank P. Kuchment for fruitful discussions and
numerous helpful comments.

This work was partially supported by the NSF/DMS grant NSF-0312292 and by the
DOE grant DE-FG02-03ER25577.

\end{document}